# Tileset As: 3 squares with local rules for non-periodic tiling

Vincent Van Dongen

Independent researcher, Canada; vincent.vandongen@gmail.com

## Abstract

This paper presents a tileset of 3 squares with local constraints on their borders and corners that enforce non-periodic tiling. We start with a description of the tileset and we demonstrate that it can tile the entire plane non-periodically creating interesting patterns. Rules are also proposed to generate the tiling. They make use of 9 supertiles with border color constraints only.

## Introduction

Wang tiles are square dominoes that can't rotate or reflect [1]. Ten years ago, an aperiodic set of 11 tiles was published [2], proven to be the one of minimal size. In this document, we explore a new set of square tiles with local constraints. These are not Wang tiles as they can rotate and both borders and corners impose constraints on their neighbors. This work is in continuation to the one done on tileset Ax [3]. It also relates to previous work on square tiles with marks for local constraints such as [4][5].

## A tileset of 3 squares with local constraints

Tileset As is a new set of tiles that consists of 3 squares, here blue, red and green. Each border is identified by a color  that takes two values (i.e. 0 or 1, or in this paper, beige or brown) and a value referred to as its thickness. Each tile also contains values at its corners. These values are constant on each tile as they are part of their specification, used by the rules to position them next to each other.

A tiler is allowed to rotate the tiles to position them in the desired orientation but rules must be followed:

- Border colors of adjacent tiles must be the same (either beige or brown);
- Sum of adjacent border thickness must be either 6 or 8;
- Sum of the 4 internal corner values must be 4, 6 or 9.

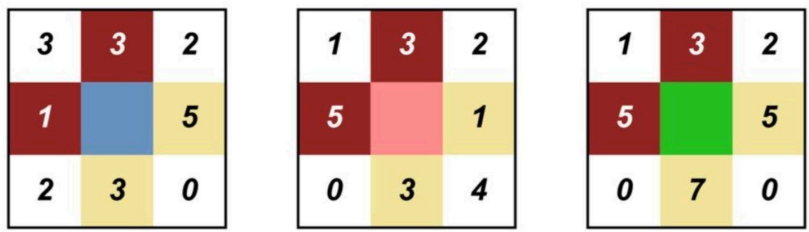

**Figure 1:** *The 3 square tiles of As with colors and values on borders and corners for local constraints.*

Figure 2 gives an example of a As tiling of size 17x17, i.e. 289 square tiles.

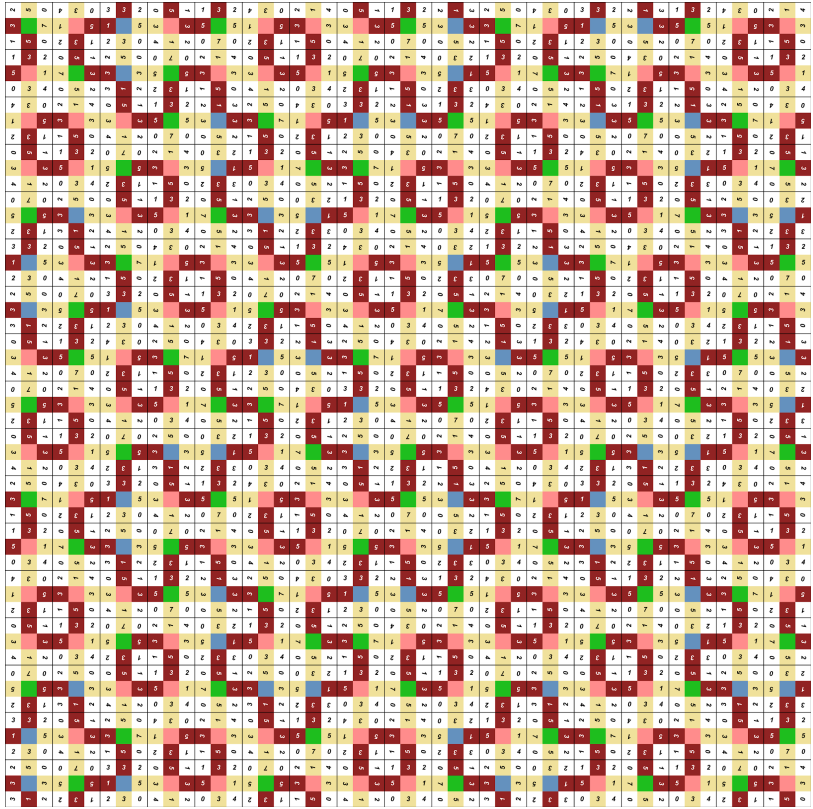

*Figure 2: Example of a 17x17 tiling As.*

## Exhaustive list of all As tilings of size 1x2

When 2 tiles are adjacent to each other, they need to follow two rules: the border color constraint and the border thickness constraint. The following figure lists all possible arrangements of 2 tiles. With all such tilings of size 2x1, one can cover the entire plane.

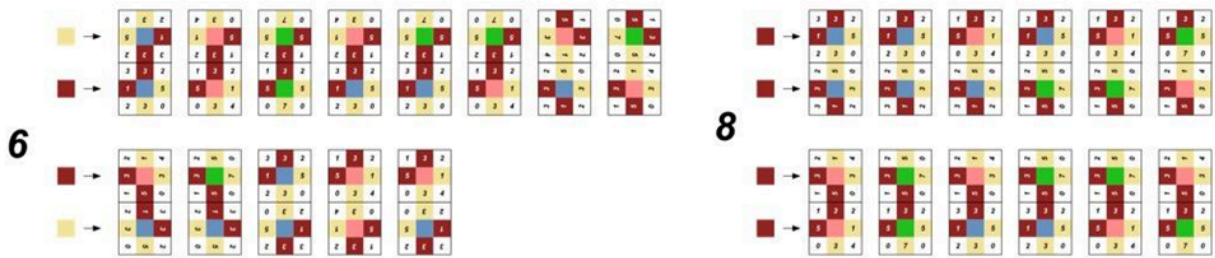

*Figure 3: All 1x2 tilings. On the left, the ones of adjacent border thickness 6. On the right, the ones of thickness 8. These two categories are incompatible for 4x4 tilings as their border color can't match.*

What appears is that, when two lines have the same border color value, their border thickness is 8, else their border thickness must be 6. This means that for a tiling to be valid, the border_thickness must be of constant value along any line of the grid.

This implies that As tilings can be structured into rectangular zones delimited by a constant border_thickness. See here below, for a subset of the above example, the zones delimited by border_thickness of value 6. Inside each zone, border_thickness is always of value 8.

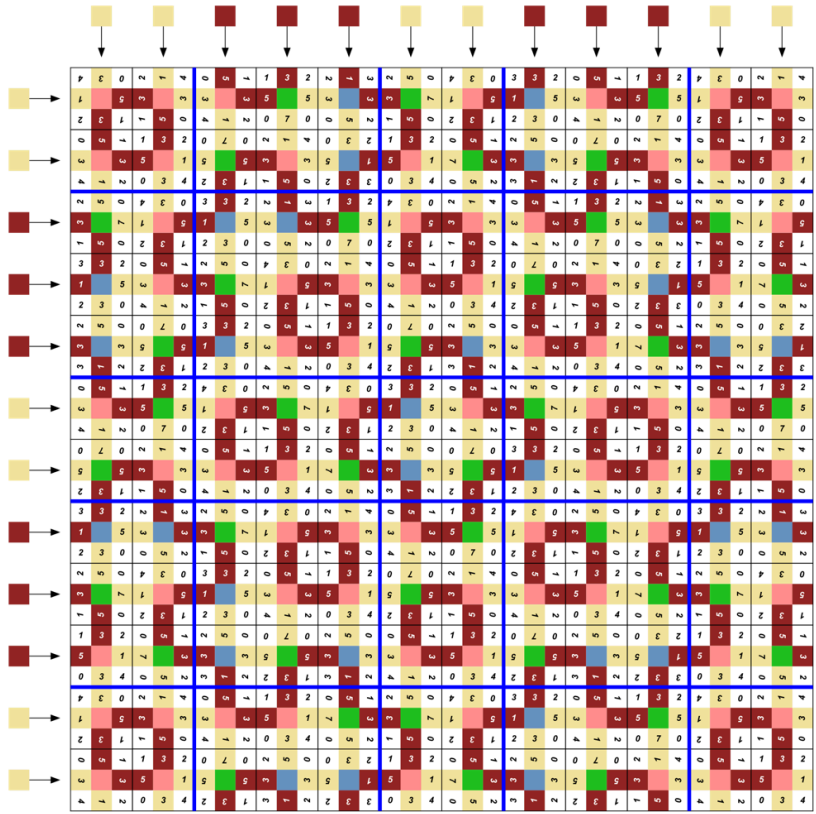

***Figure 4:*** *Along the blue lines, all border thickness is 6. This defines rectangular zones.*

## Tilings delimited by constant border thickness (or supertiles)

We searched for all the zones of As delimited by border thickness of value 6. See Figure 5. One can find by brute force that no other zone exists that could lead to valid tilings (i.e that could be further surrounded by zones). Because this set can be used to tile the entire plane, we call them supertiles.

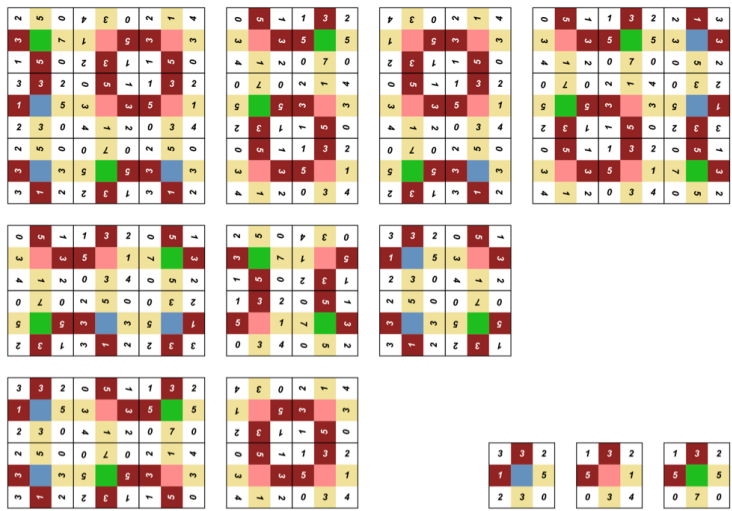

***Figure 5:*** *All the zones/supertiles of border thickness 6 that can be surrounded by zones/supertiles.*

## A set of 9 dominoes equivalent to As supertiles

By construction, the rules on the thickness of adjacent supertiles is now always 6. Also, the corners constraint on the borders of the supertile is always 6 as well. And we can show that when these two constraints are met, automatically the sum of the corners of the supertile is always 9, hence always valid. This is demonstrated by brute force. There exists no tiling of size 2x2 that satisfies the following conditions: color condition is met, border thickness is 6, corner value equals 6 (instead of 9). This means that the corners value can be removed from the specification of the tile. We are now left with a set of nine colored dominoes (with no corner constraint) that we can show can tile the entire plane just like As does (i.e. with fixed points). See figure below.

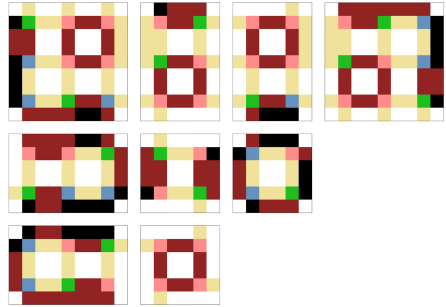

**Figure 6:** *A set of 9 dominoes (with no corner constraint) that is equivalent to As supertiles.*

## As can cover the entire plane non-periodically

Rules were developed to tile with these supertiles. In fact, each rule transforms one supertile with a large tile of the same geometric property (i.e. fixed points). This means that the rules can be used iteratively and tileset As can cover the entire plane. See Figure 7.

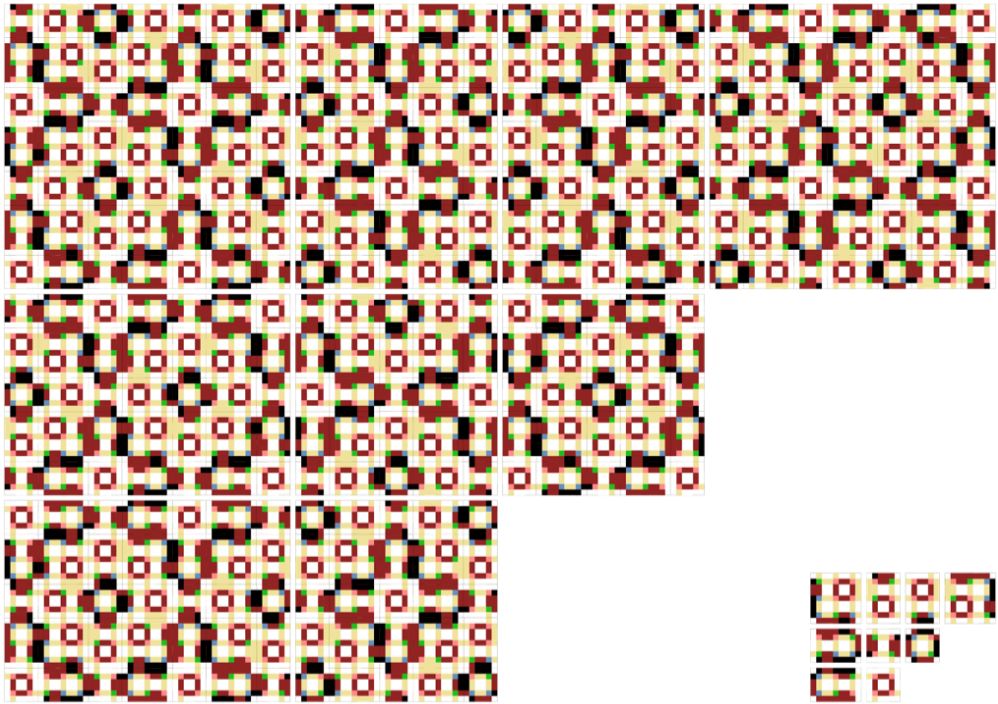

**Figure 7:** *The 9 dominoes and their corresponding fixed-point rules.*

The demonstration that tileset As has a non-periodic solution on the entire plane makes use of the bicolor grid on the 9 supertiles and fixed points rules. We assign one color for Border_thickness=6 (the boundaries of the supertiles) and one color for Border_thickness=8 (lines in blue here below).

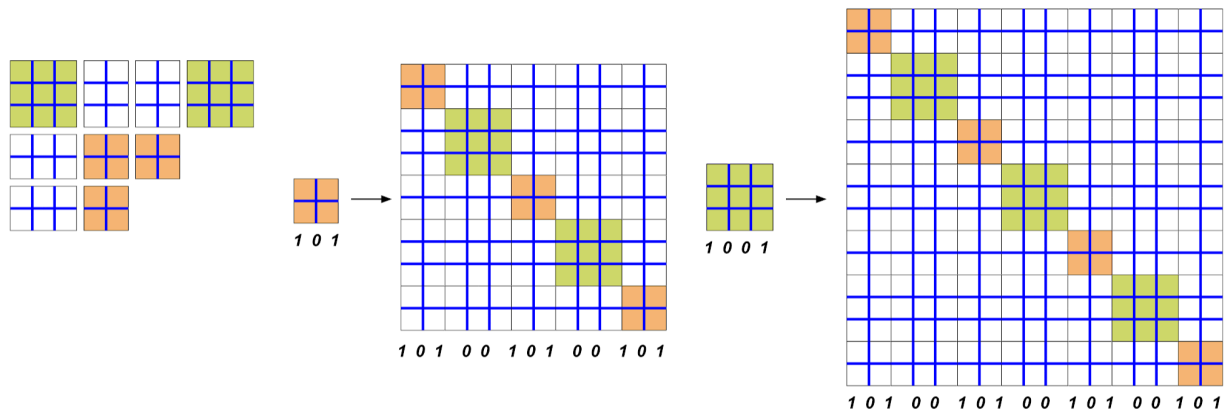

***Figure 8:*** *Example of how a supertile keeps its geometric properties after applying its rule.*

The transformation of 2x2 squares by their fixed point rule is shown, as well as the one of 3x3 squares. Blue lines are either single 0 between delimiter 1's (i.e. 1 0 1) or in pair (00) between 1's. At each iteration, 0 becomes 0 00 0 00 0, and in parallel 00 becomes 0 00 0 00 0 00 0. The sequence that it creates is well-known OEIS A159684 that is non-periodic.

## Conclusion

Tileset As is a set of 3 squares with colors and values on corners and borders to force local tiling constraints. We demonstrate in this paper that the tiling can cover the entire plane in a non-periodic manner. For this, we derived a set of 9 colored dominoes directly from As and found fixed points for these supertiles, creating non-periodic tilings this way. A brute force approach, not presented here, was able to prove that tiling with only 8 of them could not tile the plane. The proof that As is aperiodic is on its way.

## Acknowledgement

Special thanks to Pierre Gradit and Jean-François Husson for their precious comments and great support.